\documentclass[twoside,12pt]{article}
\usepackage{amsmath,amsthm,amssymb,amscd}
\title {Kernel of the variation operator\\
and periodicity of open books}
\author{Louis H. Kauffman and Nikolai A. Krylov}

\begin{document}

\date {}

\newtheorem{thm}{Theorem}
\newtheorem{cor}{Corollary}
\newtheorem{lem}{Lemma}
\newtheorem{claim}{Claim}
\newtheorem{dfn}{Definition}
\newtheorem{prop}{Proposition}

\def\Dbl            {{\cal D}}
\def\DblX           {{\cal D_A}X}
\def\DX             {{\cal D}X}
\def\DblF           {{\cal D}F}
\def\DblFi          {{\cal D}\hat{F_i}}
\def\H              {{\cal H}}
\def\V(f)           {{\rm Var(f)}}
\def\vf             {\varphi}
\def\Va             {{\rm Var}}
\def\Im(var)        {{{\rm Hom}(H_n(F,\partial F),~H_n(F))}}
\def\K(V)           {\tilde{\pi}_0V {{\rm Diff}(F,\partial)}}
\def\mcgDF          {\tilde{\pi}_0 {{\rm Diff}}(\DblF)}
\def\mcgF           {\tilde{\pi}_0 {{\rm Diff}(F,{\rm rel}~\partial)}}
\def\G              {{G}}
\def\Rhom           {{{\rm Hom}(H_n(F,\partial F),~\G)}}
\def\Dhom           {{{\rm Hom}(H_n(\DblF),~\G)}}
\def\Fhom           {{{\rm Hom}(H_n(F_+),~\G)}}
\def\kerDF          {\tilde{\pi}_0S {{\rm Diff}(\DblF)} }
\def\nat            {\mathbb N}
\def\int            {\mathbb Z}
\def\CCF            {\bar{\cal F}}
\def\RMCGD          {\tilde{\pi}_0 {\rm Diff}(DM)}
\def\ps-is          {\stackrel{\partial}{\simeq}}
\def\rmcgF          {\tilde{\pi}_0 {{\rm Diff}(F,~{\rm rel}~K)}}
\def\mcgM           {\tilde{\pi}_0 {{\rm Diff}(M)}}
\def\kerM           {\tilde{\pi}_0S {{\rm Diff}(M)} }
\def\zint           {\mathbb Z}
\def\cn             {\mathbb C}
\def\da             {\downarrow}
\def\lt             {\lefteqn}
\def\lra            {\longrightarrow}
\def\ra             {\rightarrow}
\def\hra            {\hookrightarrow}
\def\lmt            {\longmapsto}
\def\Aut            {{\rm Aut}~H_n(M)}
\def\lam            {\lambda}
\def\del            {\delta}
\def\eps            {\epsilon}
\def\tp             {\tilde{\psi }}
\def\tj             {\tilde\jmath}
\def\tf             {\tilde{\vf }}

\maketitle

\parskip=2mm

\begin{abstract}

{We consider a parallelizable $2n$-manifold $F$ which has the
homotopy type of the wedge product of $n$-spheres and show that
the group of pseudo-isotopy classes of orientation preserving
diffeomorphisms that keep the boundary $\partial F$ pointwise
fixed and induce the trivial variation operator is a central
extension of the group of all homotopy $(2n+1)$-spheres by
$H_n\bigl(F; S\pi_n(SO(n))\bigr)$. Then we apply this result to
study the periodicity properties of branched cyclic covers of
manifolds with simple open book decompositions and extend the
previous results of Durfee, Kauffman and Stevens to dimensions 7
and 15.}
\end{abstract}


\section{Introduction and the Results}

An open book decomposition of a manifold $M^{m+1}$ is a
presentation of this manifold as the union of the mapping torus
$F_{\vf}$ and the product $\partial F\times D^2$ along the
boundary $\partial F\times S^1$, where $\vf: F^m\lra F^m$ is an
orientation preserving diffeomorphism which fixes the boundary
$\partial F$ pointwise. Open book structures have been used in the
study of various topological problems (for short historical
overviews see \S2 of \cite{Quinn} or Appendix by E. Winkelnkemper
in \cite{Ranicki}), and in particular in the study of the isolated
complex hypersurface singularities. Let $f:~(\cn^{n+1},0)\lra
(\cn,0)$ be a polynomial mapping with the only singular point at
the origin and with zero locus $V=\{z\in\cn^{n+1}|f(z)=0\}$.
Consider the intersection of $V$ with a small sphere centered at
the origin $K:=V\cap S_{\eps}^{2n+1}$. J. Milnor has shown in
\cite{Mil2} that the mapping
$$
\Phi(z):=f(z)/|f(z)|~~~~ S_{\eps}^{2n+1}\setminus K\lra S^1
$$
is the projection map of a smooth fibration such that the fiber
$F:=\Phi^{-1}(1)$ is a smooth $(n-1)$-connected parallelizable
$2n$-manifold homotopically equivalent to the wedge product of
$n$-spheres and $\partial F=K$ is $(n-2)$-connected. This gives
the open book structure to the sphere
$$
S^{2n+1}=F_{\vf}\cup (K \times D^2)
$$
Such an open book decomposition of $S^{2n+1}$ is called a simple
fibered knot and the periodicity, in $k$, of the $k$-fold cyclic
covers of $S^{2n+1}$ branched along $K$ has been studied by A.
Durfee and L. Kauffman in \cite{DK}. Later, J. Stevens (see
\cite{Stev}, Theorem 7 and Proposition 8) generalized Theorems 4.5
and 5.3 of \cite{DK} to a wider class of manifolds with simple
open book decompositions $M^{2n+1}=F_{\vf}\cup (\partial F \times
D^2)$ (an open book $M^{2n+1}$ is called {\it simple} if both $M$
and $F$ are $(n-1)$-connected and $M$ bounds a parallelizable
manifold).

{\bf Theorem I} (Stevens). {\it Let $M_k$ denote the k-fold cyclic
cover of $M^{2n+1}$ branched along $\partial F$ and $n\ne
1,3~or~7~odd$. If~ ${\rm Var}(\vf^d)=0$, then $M_k$ and $M_{k+d}$
are (orientation preserving) homeomorphic, while $M_k$ and
$M_{d-k},~k<d$ are orientation reversing homeomorphic.
Furthermore, $M_{k+d}$ is diffeomorphic to
$(\sigma_d/8)\Sigma~\#~M_k$.}

\noindent Here $\sigma_k$ is the signature of a parallelizable
manifold $N_k$ with the boundary $\partial N_k=M_k$, and $\Sigma$
is the generator of the finite cyclic group $bP_{2n+2}$ of
homotopy $(2n+1)$-spheres that bound parallelizable manifolds.
${\rm Var}(h)$ denotes the variation homomorphism of a
diffeomorphism $h: F\lra F$, which keeps the boundary $\partial F$
pointwise fixed, and defined as follows. Let $[z]\in
H_n(F,\partial F)$ be the homology class of a relative cycle $z$,
then we define ${\rm Var}(h): H_n(F,\partial F)\lra H_n(F)$ by the
formula ${\rm Var}(h)[z]:=[h(z)-z]$ (cf. \S 1 of \cite{Stev} or
\S1.1 of \cite{AGV}).

Stevens also proved topological as well as smooth periodicity for
$n$ even (see \cite{Stev}, Theorem 9):

{\bf Theorem II}. {\it If for branched cyclic covers $M_k$ of a
$(2n+1)$-manifold $M$ with simple open book decomposition ${\rm
Var}(\vf^d)=0$, then $M_k$ and $M_{k+2d}$ are homeomorphic and
$M_k$ and $M_{k+4d}$ are diffeomorphic. Moreover, if $n=2~or~6$,
then $M_k$ and $M_{k+d}$ are diffeomorphic.}

Both of the papers viewed the open book $M^{2n+1}$ as the boundary
of a $(2n+2)$-manifold and used results of C.T.C. Wall \cite{W1},
on classification of $(n-1)$-connected $(2n+1)$-manifolds. Here,
in the third section, we are dealing with the same periodicity
problems from a different point of view which is based on results
of M. Kreck \cite{Kreck} on the group of isotopy classes of
diffeomorphisms of $(n-1)$-connected almost-parallelizable
$2n$-manifolds. We give here different proofs of these two
theorems of Stevens including the cases $n=3$ and $n=7$ (see
Corollaries 2, 3, and 4 below).

As we have just mentioned, our approach is based on the results of
Kreck who has computed the group of isotopy classes of
diffeomorphisms of closed $(n-1)$-connected almost-parallelizable
$2n$-manifolds in terms of exact sequences. In the first part of
this paper we use these results to obtain a similar exact sequence
for the diffeomorphisms $f$ of a parallelizable handlebody
$F\in\H(2n,\mu,n),~n\geq 2$, that preserve the boundary $\partial
F$ pointwise and induce the trivial variation operator $\V(f) :
H_*(F,\partial F)\lra H_*(F)$. We will denote the group of
pseudo-isotopy classes of such diffeomorphisms by $\K(V) $ and
prove the following

\noindent {\bf Theorem 3.}~{\it If $n\geq 3$ then the following
sequence is exact
$$
0\lra \Theta_{2n+1} \lra \K(V) \lra {{\rm
Hom}\Bigl(H_n(F,\partial F),~S\pi_n(SO(n))\Bigr)} \lra 0
$$
If $n=2$ then $\K(V) = 0$.}

\noindent Here, by $S\pi_n(SO(n))$ we mean the image of
$\pi_n(SO(n))$ in $\pi_n(SO(n+1))$ under the natural inclusion
$SO(n)\hra SO(n+1)$ and by $\Theta_{2n+1}$ the group of all
homotopy $(2n+1)$-spheres (see \S2.2 for the details).

\noindent \underline {Remark}: Recently D. Crowley \cite{Crowl}
extended results of D. Wilkins on the classification of closed
$(n-1)$-connected $(2n+1)$-manifolds, $n=3,7$. One could use these
results together with the technique of Durfee, Kauffman and
Stevens to complete the periodicity theorems for $n=3,7$. However
our intention was to show how one can apply the higher dimensional
analogs of the mapping class group in studying this kind of
problem.

At the end we briefly mention the cyclic coverings of $S^3$
branched along the trefoil knot as an example which shows that
there is no topological periodicity in the case $n=1$.

Let $F$ be a manifold with boundary $\partial F$ and consider two
diffeomorphisms $\vf,~\psi$ of $F$ that are identities on the
boundary (in this paper we consider only orientation preserving
diffeomorphisms). As usual, two such diffeomorphisms are called
{\it pseudo-isotopic relative to the boundary} if there is a
diffeomorphism ${\cal H}:~F\times I\lra F\times I$ which satisfies
the following properties:
$$
1)~{\cal H}|_{F\times \{0\}}=\vf,~~~2)~{\cal H}|_{F\times
\{1\}}=\psi,~~~3)~{\cal H}|_{\partial F\times I}=id
$$
We will denote the group of pseudo-isotopy classes of such
diffeomorphisms by $\mcgF$. The group of pseudo-isotopy classes of
orientation preserving diffeomorphisms on a closed manifold $M$
will be denoted by $\mcgM$. There is a deep result of J. Cerf
\cite{Cerf} which allows one to replace pseudo-isotopy by isotopy
provided that the manifold is simply connected and of dimension at
least six. All our manifolds are simply connected here, so $n=2$
is the only case when we actually use pseudo-isotopy. For all
other $n\geq 3$ we will use the same notations (where tilde \~ {}
stands for ``pseudo") but mean the usual isotopy. We will call
these groups {\it the mapping class groups}.

If $M$ is embedded into $W$ as a submanifold, then the normal
bundle of $M$ in $W$ will be denoted by $\nu(M;W)$. Integer
coefficients are understood for all homology and cohomology
groups, unless otherwise stated, and symbols $\simeq$ and $\cong$
are used to denote diffeomorphism and isomorphism respectively.

\noindent {\bf Acknowledgement}: We thank the referee for helpful
comments and suggestions. The second author also would like to
express his gratitude to Professors Matthias Kreck and Anatoly
Libgober for stimulating discussions during the preparation of
this paper.

\section{Kernel of the variation operator}

\subsection{Double of a pair (X,A)}

Let $(X,A)$ be a pair of CW complexes, and consider the pair
$(X\times I,A\times I)$ (here and later $I=[0,1]$, and we denote
the boundary of $I$ by $\partial I$).

\begin{dfn}
The subspace $(X\times\partial I)\cup(A\times I)$ of $X\times I$
will be called the double of the pair $(X,A)$, and denoted by
$\DblX$.
\end{dfn}

We will denote the pair $(X\times\{0\},A\times\{0\})$ by
$(X_0,A_0)$, the product $A\times I$ by $A_+$ and the union
$(X\times\{1\})\cup A_+$ by $X_+$. Thus we can write
$\DblX=X_0\cup X_+$ and $X_0\cap X_+ = A\times\{0\}$.

\noindent \underline{Remark}: If we take the pair $(X,A)$ to be a
manifold with the boundary, then the double $\DblX$ will be the
boundary of the product $X\times I$, which is a closed manifold
with the canonically defined smooth structure (see
\cite{Munkres}). In this case we will denote the double simply by
${\cal D}X$.

Now we construct a natural homomorphism $d_*: H_*(X,A)\lra
H_*(\DblX)$. Consider the reduced suspensions of $X$ and $A$ (the
common base point is chosen outside of $X$) and the induced
isomorphism between $H_*(X,A)$ and $H_{*+1}(\Sigma X^+,\Sigma
A^+)$. The excision property induces a natural isomorphism between
$H_{*+1}(\Sigma X^+,\Sigma A^+)$ and $H_{*+1}(X\times I, \DblX)$,
and we define the homomorphism $d_*$ as the composition of these
two isomorphisms with the boundary map $\delta_{*+1}$ from the
exact sequence of the pair $(X\times I,\DblX)$:

\begin{dfn}
$$ d_q:=\delta_{q+1}\circ iso :~ H_q(X,A) \stackrel{\cong}{\lra}
H_{q+1}(X\times I, \DblX) \stackrel{\delta_{q+1}}{\lra} H_q
(\DblX) $$
\end{dfn}

The groups $H_*(X,A)$ and $H_*(\DblX,X)$ are naturally isomorphic
and we can rewrite the exact sequence of the pair $(\DblX,X)$ in
the following form:
$$
\cdots H_{q+1}(\DblX)\to H_{q+1}(X,A) \to H_q(X)
\stackrel{i_q}{\to} H_q(\DblX) \stackrel{j_q}{\to} H_q(X,A)\cdots
$$

\begin{lem}
For each $q\geq 1$ the homomorphism $d_q$ is a splitting
homomorphism of the above exact sequence and we have the following
short exact sequence that splits:
$$
0\lra H_q(X) \stackrel{i_q}{\lra} H_q(\DblX) \stackrel{j_q}{\lra}
H_q(X,A)\lra 0$$
\end{lem}
\begin{proof}
It follows rather easily from our definition of $d_q$ that for
each $q\geq 1$ the composition $j_q\circ d_q$ is the identity map
of the group $H_q(X,A)$. This property entails our lemma (cf.
\cite{RF}, chap. 5, \S 1.5).
\end{proof}

Let us consider now a homeomorphism $f: X\lra X$ which is the
identity on $A$, i.e. $f(x)=x$ for all $x\in A$. For such a map
the variation homomorphism $\V(f) : H_q(X,A)\lra H_q(X)$ is
defined for all $q\geq 1$ by the formula $\V(f) [z]:=[f(z)-z]$ for
any relative cycle $z\in H_q(X,A)$ (cf. \S 1 of \cite{Stev} or
\S1.1 of \cite{AGV}). The map $f$ also induces the map $f^{(r)}:
(X,A)\lra (X,A)$ and a map $\tilde{f}: \DblX\lra \DblX$ defined as
follows:
$$\tilde{f}(x):=\left\{
\begin{array}{rcl}
f(x)& {\rm if}& x\in X_0\\
x~& {\rm if}& x\in X_+\\
\end{array}
\right.
$$

\noindent If we denote the corresponding induced maps in homology
by $f_*,~f^{(r)}_*,~\tilde{f_*}$ then we have the following
commutative diagram:
$$
\begin{CD}
0 @>>> H_q(X) @>i_q>> H_q(\DblX) @>j_q>> H_q(X,A) @>>>0\\
@.    @VVf_*V       @VV\tilde{f_*}V     @VVf^{(r)}_*V    @.\\
0 @>>> H_q(X) @>i_q>> H_q(\DblX) @>j_q>> H_q(X,A) @>>>0
\end{CD}
$$

\begin{thm}
~\\ If~ $\V(f) = 0$, then $\tilde{f_*}$ is the identity map of
$H_q(\DblX)$ for all $q$.
\end{thm}
\begin{proof}
It follows right from the definition of $\V(f) $ that $f_* - Id =
\V(f) \circ l_*$ and $f_*^{(r)} - Id = l_*\circ\V(f) $, where
$l_*: H_*(X)\lra H_*(X,A)$ is induced by the inclusion
$(X,\emptyset) \hra (X,A)$ (cf. \S1.1 of \cite{AGV}). It is also
easy to check that the homomorphisms $\tilde{f_*}$ and $d_q$ are
connected with the variation homomorphism via the formula
$$\tilde{f_*}\circ d_q = d_q\circ Id + i_q\circ \V(f) $$
Hence if $\V(f) = 0$, then $f_*=Id,~f_*^{(r)}=Id$ and
$\tilde{f_*}\circ d_q = d_q\circ Id$. These three identities
together with $j_q\circ d_q = Id$ imply the statement.
\end{proof}

Now we restrict our attention to the case when $X$ is a smooth,
simply connected manifold of dimension at least four and
$A=\partial X$ is the boundary. Let $\varphi\in {{\rm Diff}(X,{\rm
rel}~\partial)}$ and $\tilde{\varphi} \in {{\rm Diff}(\DX)}$ be
the extension by the identity to the second half of the double.
Define the map $\omega: {{\rm Diff}(X,{\rm rel}~\partial)}\lra
{\rm Diff}(\DX)$ by the formula
$\omega(\varphi):=\tilde{\varphi}$.

\begin{thm}~\\
The map $\omega$ induces a monomorphism $\tilde{\pi}_0{{\rm Diff}
(X,{\rm rel}~\partial)}\lra \tilde{\pi}_0{\rm Diff}(\DX)$.
\end{thm}
\begin{proof}
It is easy to see that $\omega$ induces a well-defined map of
groups of pseudo-isotopy classes of diffeomorphisms, i.e., if
$\varphi'$ is pseudo-isotopic relative to the boundary to
$\varphi$ then $\omega(\varphi')$ is pseudo-isotopic to
$\omega(\varphi)$. It is obvious that for any two diffeomorphisms
$\varphi,\psi\in \rm Diff(X,rel~\partial)$,
$\omega(\varphi\cdot\psi)=\omega(\varphi)\cdot\omega(\psi)$, that
is $\omega$ induces a homomorphism which we also denote by
$\omega$.

To show that $\omega$ is actually a monomorphism we use
Proposition 1 of Kreck (see \cite{Kreck}, p. 650 for the details):
{\it Let $A^m$ be a simply-connected manifold with $m\geq 5$ and
$h\in {\rm Diff}(\partial A)$. $h$ can be extended to a
diffeomorphism on $A$ if and only if the twisted double $A\cup_h
-A$ bounds a 1-connected manifold $B$ such that all relative
homotopy groups $\pi_k(B,A)$ and $\pi_k(B,-A)$ are zero, where $A$
and $-A$ mean the two embeddings of $A$ into the twisted double.}
Suppose now that $\omega(\varphi)=\tilde{\varphi}$ is
pseudo-isotopic to the identity. Then the mapping torus
$\DX_{\tilde{\varphi}}$ is diffeomorphic to the product $\DX\times
S^1 =\partial(X\times I\times S^1)$. On the other hand we can
present $\DX_{\tilde{\varphi}}$ as the union of $X_{\varphi}$ and
$-X\times S^1$ along the boundary $\partial X\times S^1$. Since
$\partial(X\times D^2)=X\times S^1\bigcup -\partial X\times D^2$
we can paste together $X\times I\times S^1$ and $X\times D^2$
along the common sub-manifold $X\times S^1$ to obtain a new
manifold $W$, which cobounds $X_{\varphi}\bigcup-\partial X\times
D^2$. Now note that $X_{\varphi}\bigcup-\partial X\times D^2$ is
diffeomorphic to the twisted double $X\times I\bigcup_h -X\times
I$ where the diffeomorphism $h: \partial(X\times I)\rightarrow
\partial(X\times I)$ is defined by the identities:
$\left.h\right|_{X_0}=id$, and $\left.h\right|_{X_+}=\varphi$ (cf.
\cite{Kreck}, property 1) of $\tilde{W}$ on page 657). The theorem
of Seifert and Van Kampen entails that $\pi_1(W)\cong\{1\}$, and
hence $\pi_1(W,X\times I)\cong\{1\}$. To show that the other
homotopy groups are trivial it is enough to show that
$H_*(W,X)\cong \{0\}$ for all $*\ge 2$. This can be seen from the
relative Mayer-Vietoris exact sequence of pairs $(X\times I\times
S^1, X)$ and $(X\times D^2, X)$ where by $X$ we mean a fiber of
the product $X\times S^1$: $H_*(X\times
S^1,X)\stackrel{\cong}{\rightarrow} H_*(X\times I\times
S^1,X)\oplus H_*(X\times D^2,X)\rightarrow H_*(W,X)$. Thus by
Proposition 1 of \cite{Kreck}, there is a diffeomorphism of
$X\times I$ to itself that gives the required pseudo-isotopy
between $\varphi$ and $id$.
\end{proof}

\subsection{$\K(V) $ as an extension}

We now let $F\in\H(2n,\mu,n)$ be a parallelizable handlebody, that
is, a parallelizable manifold which is obtained by gluing $\mu$
$n$-handles to the $2n$-disk and rounding the corners:
$$
F=D^{2n}\cup\bigsqcup_{i=1}^{\mu} (D_i^n\times D^n)$$ We assume
here that $n\geq 2$. For the classification of handlebodies in
general, see \cite{Wall1}. Obviously $F$ has the homotopy type of
the wedge product of $n$-spheres and nonempty boundary $\partial
F$ which is $(n-2)$-connected. The Milnor fibre of an isolated
complex hypersurface singularity is an example of such a manifold.

Let us consider now $\vf\in\mcgF$ and the induced variation
homomorphism $\Va(\vf): H_n(F,\partial F)\lra H_n(F)$. This
correspondence gives a well defined map
$$\Va: \mcgF \lra \Im(var) $$
which is a derivation (1-cocycle) with respect to the natural
action of the group $\mcgF$ on $\Im(var) $ (cf. \cite{Stev}, \S2)
$$
\Va(h\circ g) = \Va(h) + h_*\circ \Va(g)$$

This formula implies that the isotopy classes of diffeomorphisms
that give trivial variation homomorphisms form a subgroup of
$\mcgF$.

\begin{dfn}
The subgroup
$$\K(V) :=\{f\in\mcgF~|~\Va(f)[z]=0,~ \forall [z]\in
H_n(F,\partial F)\}$$ will be called the kernel of the variation
operator.
\end{dfn}

In order to describe the algebraic structure of this kernel we
will use the results of Kreck \cite{Kreck} who has computed the
group of isotopy classes of diffeomorphisms of closed oriented
$(n-1)$-connected almost-parallelizable $2n$-manifolds in terms of
exact sequences. First we note that the double of our handlebody
$F$ is such a manifold.

\begin{lem}
Let $F\in\H(2n,\mu,n)$ be a parallelizable handlebody $(n\geq 2)$,
then the double $\DblF$ is a closed $(n-1)$-connected
stably-parallelizable $2n$-manifold.
\end{lem}
\begin{proof}
Since $F$ is simply connected and $\DblF=F_0\cup F_+$, we have
$\pi_1(\DblF)=0$. Then using the exact homology sequence of the
pair ($F\times I, \partial(F\times I)$) it can be easily seen that
$\DblF$ is a $(n-1)$-connected manifold. Since $F$ is
parallelizable the double will be stably-parallelizable.
\end{proof}

Next we recall the result of Kreck \cite{Kreck}. Let $M$ be a
smooth, closed, oriented $(n-1)$-connected almost-parallelizable
$2n$-manifold, $n\geq 2$. Denote by $\Aut$ the group of
automorphisms of $H_n(M,\int)$ preserving the intersection form on
$M$ and (for $n\geq3$) commuting with the function
$\alpha:~H_n(M)\lra \pi_{n-1}(SO(n))$, which is defined as
follows. Represent $x\in H_n(M)$ by an embedded sphere
$S^n\hookrightarrow M$. Then function $\alpha$ assigns to $x$ the
classifying map of the corresponding normal bundle. Any
diffeomorphism $f\in {\rm Diff}(M)$ induces a map $f_*$ which lies
in $\Aut$. This gives a homomorphism
$$
\kappa:~ \mcgM \lra \Aut,~~~[f]\longmapsto f_*
$$

The kernel of $\kappa$ is denoted by $\kerM$ and to each element
$f$ from this kernel Kreck assigns a homomorphism $H_n(M)\lra
S\pi_n(SO(n))$, where $S: \pi_n(SO(n))\lra \pi_n(SO(n+1))$ is
induced by the inclusion, in the following way. Represent $x\in
H_n(M)$ by an imbedded sphere $S^n\subset M$ and use an isotopy to
make $f|_{S^n}=Id$. The stable normal bundle $\nu(S^n)\oplus
\varepsilon^1$ of this sphere in $M$ is trivial and therefore the
differential of $f$ gives an element of $\pi_n(SO(n+1))$. It is
easy to see that this element lies in the image of $S$. This
construction leads to a well defined homomorphism (cf. Lemma 1 of
\cite{Kreck})
$$\chi:~\kerM\lra {{\rm Hom}\bigl(H_n(M),~S\pi_n(SO(n))\bigr)}$$

\noindent If $n=6$ we have $S\pi_n(SO(n))=0$, and for all other
$n\geq 3$ the groups $S\pi_n(SO(n))$ are given in the following
table (\cite{Kreck}, p. 644):

\parskip=5mm

\begin{tabular}{|c|c|c|c|c|c|c|c|c|}
\hline $n$ (mod 8) & 0 & 1 & 2 & 3 & 4 & 5 & 6 & 7\\
\hline $S\pi_n(SO(n))$ & ~$\zint_2\oplus\zint_2$~ & ~$\zint_2$~ &
~$\zint_2$
~ & ~$\zint$ ~ & ~ $\zint_2$~ & ~ 0 ~ & ~ $\zint_2$ ~ & ~$\zint$ \\
\hline
\end{tabular}

In particular, when $n\equiv 3\pmod{4}$ the homomorphism $\chi(f)$
can be defined using the Pontryagin class $p_{\frac{n+1}{4}}(M_f)$
of the mapping torus $M_f$: Take a diffeomorphism $f\in \kerM$ and
consider the projection
$$
\pi:~M_f\lra \left. M_f \right/ \{0\}\times M =\Sigma M^+
$$
It is clear from the exact sequence of Wang that the map
$i^*:~H^n(M_f)\lra H^n(M)$ is surjective (recall that $f_*=id$)
and therefore we obtain an isomorphism $\pi^*:~H^n(M)\cong
H^{n+1}(\Sigma M^+)\lra H^{n+1}(M_f)$. Next define an element
$p'(f)\in H^n(M)$ by
$p'(f):={\pi^*}^{-1}(p_{\frac{n+1}{4}}(M_f))$. It can be shown
(cf. \cite{Kerv}) that the map $f\lmt p'(f)$ is a homomorphism and
$c:=a_{\frac{n+1}{4}} (\frac{n-1}{2})!$ divides $p'(f)$, where as
always, $a_m=2$ if $m$ is odd and $a_m=1$ if $m$ is even. This
defines a map
$$
\chi':~\kerM\lra H^n(M),~~~with~~~\chi'(f):= p'(f)/c
$$

\parskip=2mm

\noindent \underline{Remark}: These two elements $\chi'(f)$ and
$\chi(f)$ belong to the isomorphic groups ${{\rm
Hom}(H_n(M),~\pi_n(SO))}~~{\rm and}~~{{\rm Hom}(H_n(M),~
S\pi_n(SO(n)))}$ respectively, and they are connected through
$\tau^*(\chi(f))=\chi'(f)$ via the homomorphism
$$
\tau^* : {{\rm Hom}\bigl(H_n(M),~S\pi_n(SO(n))\bigr)} \lra {{\rm
Hom}(H_n(M),~\pi_n(SO))}
$$
induced by the natural homomorphism $\tau: \pi_n(SO(n+1))\lra
\pi_n(SO(n+2))$. For the details the reader is referred to Lemma 2
of \cite{Kreck}.

If $M^{2n}$ bounds a parallelizable manifold and $n\geq 3$, then
Theorem 2 of \cite{Kreck} gives two short exact sequences:

\begin{equation}
\label{firstseq} 0\lra \kerM \lra \mcgM \stackrel{\kappa}{\lra}
\Aut \lra 0
\end{equation}
\begin{equation}
\label{secondseq} 0\lra \Theta_{2n+1} \stackrel{\iota}{\lra} \kerM
\stackrel{\chi}{\lra} {{\rm Hom}\bigl(H_n(M),~S\pi_n(SO(n))\bigr)} \lra 0
\end{equation}

\noindent where the map $\iota$ is induced by the identification
of each homotopy $(2n+1)$-sphere with the element of the mapping
class group $\tilde{\pi}_0 {\rm Diff}(D^{2n},{\rm rel}~\partial)$.

If $M$ is a simply connected manifold of dimension 4, Kreck has
proved that $\kappa$ is a monomorphism (\cite{Kreck}, Theorem 1).

Let $F\in\H(2n,\mu,n)$ be a parallelizable handlebody as above,
and $\DblF$ be the corresponding double. First assume that $n=2$
and $\vf\in \K(V) $, then it follows from our Theorems 1 and 2 and
Theorem 1 of Kreck \cite{Kreck} that $\tilde{\vf}$ is the trivial
element of $\mcgDF$, and therefore $\vf$ is the identity of
$\mcgF$.

\noindent \underline{Remark}: In this case, the handlebody $F$
doesn't have to be parallelizable and the kernel of the variation
operator $\K(V) $ will be trivial for any simply connected
4-manifold $F$.

\noindent Next we consider the case when $n\geq 3$ and denote the
group $S\pi_n(SO(n))$ by $\G$. Recall also that we can assume that
$\DblF=F\cup F_+$. Since $F$ is $(n-1)$-connected and the boundary
$\partial F$ is $(n-2)$-connected, the universal coefficient
theorem together with the cohomology exact sequence of the pair
$(\DblF,F_+)$ and the excision property give us the following
short exact sequence:
\begin{equation}
\label{thirdseq} 0\to \Rhom \stackrel{j^*}{\to} \Dhom
\stackrel{i^*}{\to} \Fhom \to 0
\end{equation}
where $i:F_+\hra \DblF$, $j: (\DblF,\emptyset)\hra (\DblF,F)$ are
inclusions and $i^*$ and $j^*$ are the corresponding induced maps.

\begin{lem}
$i^*(\chi(\tilde{\vf}))$ is the trivial map for any $\vf\in\mcgF$.
\end{lem}
\begin{proof}
Take any $[z]\in H_n(F_+)$, then we have
$i^*(\chi(\tilde{\vf}))[z]=\chi(\tilde{\vf})[i_*(z)]$. Since
$H_n(F)\cong\pi_n(F)$ we can present our $n$-cycle $[z]$ by an
imbedded $S^n\hra F_+$ and we can also assume that the normal
bundle of such a sphere is contained in $F_+$. We have defined
$\tilde{\vf}$ as the identity on $F_+$ and this implies
$\chi(\tilde{\vf})[i_*(z)]=0$ as required.
\end{proof}

\noindent Now we define a homomorphism $\chi_r: \K(V) \to \Rhom$.
Take any $\vf\in \K(V) $ then $\tilde{\vf}\in \kerDF$ (recall
Theorem 1 above) and $\chi(\tilde{\vf})\in \Dhom$. Since
$i^*(\chi(\tilde{\vf}))=0$ there exists unique $h\in \Rhom$ such
that $j^*(h)=\chi(\tilde{\vf})$.

\begin{dfn}
We define the map $\chi_r: \K(V) \to \Rhom $ by the formula
$\chi_r(\vf):=h$.
\end{dfn}
\noindent It is clear that $\chi_r$ is
a homomorphism. Here we also consider the map $\iota_r:
\Theta_{2n+1}\lra \mcgF$ defined as in (\ref{secondseq}) above:
present any homotopy $(2n+1)$-sphere $\Sigma'$ as the union of two
disks via a diffeomorphism $\psi\in \tilde{\pi}_0 {\rm
Diff}(S^{2n})\cong \tilde{\pi}_0 {\rm Diff}(D^{2n},{\rm
rel}~\partial ) \cong \Theta_{2n+1}$ then take a disk $D^{2n}$
embedded into ${\rm int}(F)$ and define the diffeomorphism of $F$
by the formula
$$
\iota_r(\Sigma')(x):=\left\{
\begin{array}{rcl}
\psi(x)& {\rm if}& x\in D^{2n}\hra F\\
x~& ~ & {\rm otherwise} \\
\end{array}
\right.
$$

\noindent It is obvious that ${\rm Im(\iota_r)}\subset \K(V) $. Now we describe $\K(V)
$ as a central extension of the group $\Theta_{2n+1}$ by $H^n(F,\partial F;\G)\cong
H_n(F;\G)$.

\begin{thm} ~\\
If $n=2$ then $\K(V) = 0$, and for all $n\geq 3$ the following sequence is exact
\begin{equation}
\label{kervar} 0\lra \Theta_{2n+1} \stackrel{\iota_r}{\lra} \K(V)
\stackrel{\chi_r}{\lra} \Rhom \lra 0
\end{equation}
\end{thm}
\begin{proof}
We have mentioned already that if $n=2$, the kernel of the variation operator is
trivial. Assume now that $n\geq 3$. It follows from Theorems 1 and 2 above that the
inclusion map $\omega: {\rm Diff}(F,{\rm rel}~\partial)\to {\rm Diff} (\DblF)$ induces
a monomorphism $s_{\omega}: \K(V) \to \kerDF$. Since the composition $s_{\omega}\cdot
\iota_r$ coincides with the injective map $\iota$ from the exact sequence
(\ref{secondseq}), we see that our $\iota_r$ is injective too. It is also clear that
$\rm Im(\iota_r)\subset Ker(\chi_r)$. Consider now any $\vf\in \rm Ker(\chi_r)$, then
$\chi( s_{\omega}(\vf)) = j^*(\chi_r(\vf)) = 0$, where $j^*$ is as in
(\ref{thirdseq}). Thus $s_{\omega}(\vf) \in \rm Ker(\chi)\cong\Theta_{2n+1}\cong
Im(\iota)$ and since $s_{\omega}$ is a monomorphism we have $\vf\in \rm Im(\iota_r)$
as required.

To prove that $\chi_r$ is an epimorphism it is enough to show that for a set of
generators $\{g_1,\ldots, g_m\}$ of $\Rhom$ the group $\K(V) $ contains
diffeomorphisms $\{\vf_1,\ldots, \vf_m\}$ such that $\chi_r(\vf_j)=g_j, ~
j\in\{1,\ldots,m\}$. Recall that $F=D^{2n}\cup\bigsqcup_{i=1}^{\mu} (D_i^n\times
D^n)$ and $H_n(F,\partial F)\cong \zint^{\mu}$. We can choose the following embedded
disks $d_i\hra F,~ i\in\{1,\ldots,\mu\}$, as a basis of this homology group:
$$
d_i:= \{0\}_i\times D^n\hra D^n_i\times D^n \hra F$$ (here $\{0\}_i$ is the center of
the $i^{th}$ handle core disk $D^n_i$). Take a generator $x$ of $G$ and consider the
homomorphism $ g_{xi}: H_n(F,\partial F)\lra G$ defined by the formula
$$
g_{xi}[d_k]:=\left\{
\begin{array}{rcl}
x & {\rm if} & k=i\\
0 & {\rm if} & k\ne i
\end{array}
\right. k\in\{1,\ldots,\mu\}
$$ end extended linearly to the whole group. The set of such homomorphisms
obviously generates $\Rhom$. Now we will use an analog of the Dehn twist in higher
dimensions to construct the diffeomorphism $\vf_{xi}$ (cf. \cite{Wall1}, Lemma 12).

For each disk $d_k$ consider the ``half-handle" $=(\frac{1}{2}D_k^n)\times D^n$ and
notice that the closure of the complement to all these ``half-handles" in $F$
$$
\CCF:=cl(F\setminus\bigsqcup_{k=1}^{\mu}(\frac{1}{2}D_k^n)\times D^n)$$ is
diffeomorphic to the closed $2n$-disk $D^{2n}$, and the intersection of each
``half-handle" with the boundary $\partial \CCF\simeq S^{2n-1}$ is $\partial
(\frac{1}{2}D_k^n)\times D^n\simeq S_k^{n-1}\times D^n$. We take a smooth map $\vf_x:
(D^n,S^{n-1})\lra (SO(n),id)$ that sends a neighborhood of $S^{n-1}$ to $id$ and
represents an element $[\vf_x]\in \pi_n(SO(n))$ such that $S([\vf_x])=x$ and define
the diffeomorphism $\vf_{xi}|_{\bigsqcup_{k=1}^{\mu}(\frac{1}{2}D_k^n)\times D^n}$ by
the formula
\begin{equation}
\label{twist}
\vf_{xi}(t,s):=\left\{
\begin{array}{ccl}
(\vf_x(s)\circ t,s) & {\rm if} & (t,s)\in (\frac{1}{2}D_i^n)\times
D^n\\
(t,s) & {\rm if} & (t,s)\in (\frac{1}{2}D_k^n)\times D^n~and~k\ne
i
\end{array}
\right.
\end{equation} In particular, this gives a diffeomorphism $\phi\in
{\rm Diff(}\partial \CCF)$ which is defined on $S_i^{n-1}\times D^n\hra \partial \CCF$
by restricting $t$ to the boundary of $\frac{1}{2}D_i^n$ (see (\ref{twist}) above)
and by the identity everywhere else. We will show now that $\phi$ is isotopic to the
identity. Consider the handlebody
$$
F_i:=D^{2n}\cup (D_i^n\times D^n)=cl(F\setminus\bigsqcup_{k=1,k\ne
i }^{\mu}(\frac{1}{2}D_k^n)\times D^n)
$$
and denote by $\hat{F_i}$ the manifold obtained from $F_i$ by removing the open disk
$\frac{1}{2}D^{2n}$ from $D^{2n}$. Hence $\partial \hat{F_i} \simeq \partial
F_i\sqcup S^{2n-1}$. The first equation of (\ref{twist}) together with the identity
map define a diffeomorphism $\Phi$ of $\hat{F_i}$ such that $\Phi|_{S^{2n-1}} = \phi$
and $\Phi|_{\partial F_i} = Id$. We use the identity again to extend this $\Phi$ to a
diffeomorphism $\tilde{\Phi}$ of $\DblFi$ where
$$
\DblFi:=\DblF_i\setminus\frac{1}{2}D^{2n}\simeq \hat{F_i}\bigcup_{\partial F_i}
F_i~~~~~and~~~~\tilde{\Phi}|_{F_i}=Id,~\tilde{\Phi}|_{\hat{F_i}}=\Phi
$$

Thus $\phi$ is the restriction of $\tilde{\Phi}$ to the boundary $\partial
\DblFi=S^{2n-1}$ and hence can be considered as an element of the inertia group of
${\cal D}F_i$ (cf. \cite{Kreck}, Proposition 3). Now it follows from Lemma 2 above
and results of Kosinski (\cite{Kos}, see \S3) and Wall \cite{W0} that $\phi$ is
isotopic to the identity. In particular, we can use this isotopy on $S^{2n-1}\times
[\frac{1}{2},\frac{1}{4}]\subset \frac{1}{2}D^{2n}$ to extend the diffeomorphism
$\vf_{xi}|_{\bigsqcup_{k=1}^{\mu}(\frac{1}{2}D_k^n)\times D^n}$ to a diffeomorphism
of the whole handlebody $F$. Denote the result of this extension by $\vf_{xi}$.
Clearly $\vf_{xi}\in \mcgF$, and we leave it to the reader to check that
$\chi_r(\vf_{xi})=g_{xi}$.
\end{proof}

\begin{cor}We have the following commutative diagram
$$
\begin{array}{ccccccccc}
&&&&0&&0&&\\
&&&&\da&&\da&&\\
0 &\ra&\Theta_{2n+1}&\stackrel{\iota_r}{\to}&\K(V)
&\stackrel{\chi_r}{\lra} &\Rhom &\ra &0\\
&& \updownarrow\lt{\equiv} &&\da\lt{s_{\omega}} &&\da\lt{j^*} &\\
0 &\ra &\Theta_{2n+1} &\stackrel{\iota}{\to} &\kerDF &\stackrel{\chi}{\lra}
& \Dhom &\ra &0\\ &&&&\da\lt{i^*\cdot\chi} &&\da\lt{i^*} &\\
&&&&\Fhom &\stackrel{\equiv}{\longleftrightarrow} &\Fhom &\\
&&&&\da&&\da&&\\
&&&&0&&0&&
\end{array}
$$ where all horizontal and vertical sequences are exact.
\end{cor}
\begin{proof} The standard diagram chasing procedure is left to the reader.
\end{proof}

\noindent {\tt Example:} Consider the case when $F=S^3\times D^3$. Then
$\DblF=S^3\times S^3$, $\Rhom\cong G\cong \zint$, $\Theta_{7}\cong\zint_{28}$ and
$\kerDF\cong{\cal H}_{28}$, that is the factor group of the group $\cal H$ (upper
unitriangular $3\times 3$ matrices with integer coefficients) modulo the cyclic
subgroup $28\zint$, where $\zint$ is the center of $\cal H$ (cf. \cite{Fried} or
\S1.3 of \cite{Kryl}). Thus $\K(V) \cong S\pi_3(SO(3))\oplus\Theta_7\cong
\zint\oplus\zint_{28}$ and the first vertical short exact sequence from the previous
corollary can be written as follows
$$
0\lra S\pi_3(SO(3))\oplus\Theta_7 \lra \tilde{\pi}_0S{{\rm Diff}(S^3\times S^3)} \lra
S\pi_3(SO(3)) \lra 0.
$$
Such exact sequence was obtained by J. Levine (\cite{Levine}, Theorems 2.4 and 3.3)
and H. Sato (\cite{Sato}, Theorem II) for the group $\tilde{\pi}_0S{{\rm
Diff}(S^p\times S^p)}$.

\section{Manifolds with open book decompositions}

\subsection{Periodicity in higher dimensions}

In this section we will apply our exact sequence (\ref{kervar}) to
study the periodicity of branched cyclic covers of manifolds with
open book decompositions.

\begin{dfn}
We will say that a smooth closed $(m+1)$-dimensional manifold $M$
has an open book decomposition if it is diffeomorphic to the union
$$
M\simeq F_{\vf}\cup_r (\partial F\times D^2)$$ where $F$ is
m-dimensional manifold with boundary $\partial F$, $\vf \in {{\rm
Diff}(F,{\rm rel}~\partial)}$ is an orientation preserving
diffeomorphism of $F$ that keeps the boundary pointwise fixed,
$F_{\vf}$ is the mapping torus of $\vf$
$$
F_{\vf}:= \left. F\times[0,1]\right/ (x,0)\sim (\vf(x),1)
$$
and $r:\partial F_{\vf}\lra \partial F\times S^1$ is a
diffeomorphism that makes the following diagram commute
$$\begin{CD}
\partial F_{\vf} @>in>> F_{\vf}\\
@VVrV  @VV\pi V\\
\partial F\times S^1 @>p_2>> S^1
\end{CD}
$$ (here $p_2$ is the projection onto the second factor and
$\pi$ is the bundle projection of the mapping torus onto the base
circle).
\end{dfn}
Such a union is also called the relative mapping torus with page
$F$ and binding $\partial F$ (cf. \cite{Quinn} or \cite{Stev}).
When $M$ has dimension $(2n+1)$ and $F$ has the homotopy type of a
$n$-dimensional CW-complex, it is said that the page is {\it
almost canonical}. The diffeomorphism $\vf$ is called the
geometric monodromy and the induced map $\vf_*:~H_n(F)\ra H_n(F)$
is the (algebraic) monodromy. If instead of $\vf$ we take some
positive power of this diffeomorphism, say $\vf^k$, we obtain the
$k$-fold cyclic cover $M_k$ of $M$, branched along $\partial F$,
i.e.
$$
M_k=F_{\vf^k}\cup_r (\partial F \times D^2)
$$

It was shown in \cite{DK} (Theorem 4.5) that if a fibered knot
$\partial F$ is a rational homology sphere and $\vf^d=id$ for some
$d>0$, then the $k$-fold cyclic covers $M_k$ of $S^{2n+1}$
branched along $\partial F$ have the periodic behavior in $d$. In
case of the links of isolated complex polynomial singularities
these restrictions on $\partial F$ and $\vf$ are equivalent to the
condition ${\rm Var}(\vf^d)=0$.

\noindent \underline{Remarks}:\\
i) Notice that the conditions $\vf^d_*=id$ and $\partial F$ is a
rational homology sphere imply that ${\rm Var}(\vf^d)=0$, but the
converse is not true (see \cite{Stev},
p. 231).\\
ii) Proposition 3.3 of \cite{Kauffman} proves that an open book
$M^{2n+1}$ with page $F$ and monodromy $\vf$ is a homotopy sphere
if and only if ${\rm Var}(\vf)$ is an isomorphism.

In addition to the almost canonical page requirement we will need
to assume more about $M$ (cf. \cite{Stev}, \S3 p.232), i.e. we
assume from now on that $M$ has a {\it simple open book
decomposition}. It implies, in particular, that $M$ bounds a
simply connected parallelizable manifold. We will also assume that
$n\geq 3$, $\pi_1(\partial F)=1$ and ${\rm Var}(\vf^d)=0$ for some
$d\ge 1$ (where $\vf$ is the diffeomorphism that gives $M$ the
open book structure). A parallelizable simply connected manifold
bounded by $M$ will be denoted by $N$.

Before we give proofs of the periodicity theorems (Corollaries 2,
3 and 4 below) we will first obtain some auxiliary results. It is
clear that $F$ is also a parallelizable manifold. Take now any
$z\in H_n(F,\partial F)\cong \pi_n(F,\partial F)$ and choose an
embedded disk $(D^n,\partial D^n)\hra (F,\partial F)$ that
represents this class $z$. Inside of $\DblF=F_0\cup F_+$ we
consider the double $\Dbl D^n= D^n_0\cup D^n_+\simeq S^n\hra
\DblF$, and since the boundary $\partial D^n=S^{n-1}\subset
\partial F$ has trivial normal bundle in $\partial F$ we can add
to $F_0$ one $n$-handle along this sphere to obtain the manifold
$F_0(z):= F_0\cup (D^n_+\times d^n)$. As we have done above, we
can extend a diffeomorphism $\vf\in {{\rm Diff}(F,{\rm
rel}~\partial)}$ to a diffeomorphism $\vf_z\in {{\rm
Diff}(}F_0(z),{\rm rel}~\partial)$ using the identity on
$D^n_+\times d^n$. Then we obviously have $\Dbl D^n\hra F_0(z)\hra
\DblF$ and $\vf_z = \tilde{\vf}|_{F_0(z)}$.

\begin{lem}
The mapping torus $F_{\vf}$ of $\vf$ is framed if and only if the
mapping torus ${F_0(z)}_{\vf_z}$ of $\vf_z$ is framed.
\end{lem}
\begin{proof}
We will show that any framing of $F_{\vf}$ can be extended to a
framing of ${F_0(z)}_{\vf_z}$. The other direction is trivial.
Since $\vf$ is the identity on the boundary, we have
$S^{n-1}\times S^1\hra \partial F_{\vf} = (\partial F)\times S^1$,
where $S^{n-1}$ is the boundary of our relative homology class
$z$. We can assume that $F$ has a collar $\partial F\times [0,1]$
and $\vf$ is the identity map on this collar. Now we have $D^n\hra
F\hra F\cup (\partial F\times [0,1])$ and we use the disk theorem
to change $\vf$ by an isotopy to a diffeomorphism $\vf'$ such that
$\vf'|_{D^{2n}}=\vf'|_{\partial F\times 1} =id$ and $D^n\subset
{\rm int}(D^{2n})\subset F\cup(\partial F\times [0,\frac{1}{2}])$.
Then clearly $F_{\vf'}\simeq F_{\vf}$ and $D^n\times S^1\hra
F_{\vf'}$ with the trivial normal bundle. Furthermore since
$S^{n-1}\times [0,1]\hra \partial F\times [0,1]$ with trivial
normal bundle too, we can connect $\partial D^n\times S^1 \hra
F_{\vf'}$ with $S^{n-1}\times S^1 \hra (\partial F\times 1)\times
S^1=\partial (F_{\vf'})$, using the collar $(S^{n-1}\times
[0,1])\times S^1$. This implies that the trivial normal bundle of
$S^{n-1}\times S^1$ in $\partial (F_{\vf})$ comes from the trivial
normal bundle of $D^n\times S^1$ embedded into $F_{\vf}$. Now
notice that the mapping torus ${F_0(z)}_{\vf_z}$ is the union of
$F_{\vf}$ and $D_+^n\times d^n\times S^1$ along $S^{n-1}\times d^n
\times S^1 \hra \partial (F_{\vf})$. Therefore the restriction of
the framing of $F_{\vf}$ to $S^{n-1}\times d^n \times S^1 =
(\partial D^n)\times d^n \times S^1$ (where $D^n\times d^n \times
S^1\hra F_{\vf}$) can be extended to a framing of
${F_0(z)}_{\vf_z}$.
\end{proof}

\begin{thm}
(n is odd, $\ne 1$)\\
Suppose $[\psi]\in \K(V) $ and $M^{2n+1}\simeq F_{\psi}\cup_r
(\partial F\times D^2)$ bounds a parallelizable manifold $N$. Then
$\chi_r(\psi) = 0$.
\end{thm}
\begin{proof}
It is enough to show that $\chi_r(\psi)[z] = 0$ for an arbitrary
relative homology class $z\in H_n(F,\partial F)$. As we just did
above, we represent such a class by an embedded disk
$(D^n,\partial D^n) \hra (F,\partial F)$ and take the double $\Dbl
D^n = S^n\hra \DblF$. We will denote this double by $dz$ ({\it to
avoid cumbersome notations we denote by $dz$ both the homology
class and the embedded sphere $\Dbl D^n$ that represents this
class}) and its normal bundle in $\DblF$ by $\nu(dz;\DblF)$
respectively. Note that $\nu(dz;\DblF)$ is trivial. The proof now
splits into two parts.\\
1)~ $n\geq 5$: It is clear that ${\psi_z}_* = id$ on $H_n(F_0(z))$
and we can isotope $\psi_z$ to a diffeomorphism $\psi'_z$ such
that ${\psi'_z}|_{dz} = id$ (see \cite{Haef}). Extending this new
diffeomorphism by the identity to the diffeomorphism
${\tilde\psi}'\in {{\rm Diff}(\DblF)}$ we obtain an element of
$\kerDF$ which pointwise fixes $dz$ and maps $F_0(z)$ to itself.
Now it follows from the commutative diagram of Corollary 1 that it
is enough to show that $\chi({\tilde \psi}')[dz]=0$. Since by
Lemma 4 the mapping torus ${F_0(z)}_{\psi'_z}$ is framed, the
normal bundle $\nu(dz\times S^1;\DblF_{{\tilde\psi}'})$ is stably
trivial. Since $n$ is odd, the map $\G = S\pi_n(SO(n)) \hra
\pi_n(SO(n+1))\lra \pi_n(SO(n+2))$ is a monomorphism (see
\cite{Wall3}) and therefore the map
$$
l^*: {{\rm Hom}(H_n(\DblF),~\G)} \lra {{\rm
Hom}(H_n(\DblF),~\pi_n(SO))}$$ is a monomorphism too. Hence
$l^*(\chi({\tilde \psi}'))[dz]$ is the obstruction to triviality
of the stable normal bundle $\nu(dz\times
S^1;\DblF_{{\tilde\psi}'})$ and since this bundle is trivial we
have $\chi({\tilde \psi}')[dz] =0$, as required.\\ 2)~ $n\equiv
3\pmod{8}$: Since $\partial N = M^{2n+1}\simeq F_{\psi}\cup_r
(\partial F\times D^2)$ and $\partial(F\times D^2) = (\partial F
\times D^2)\cup (F\times S^1)$, we can paste together the
manifolds $N$ and $F\times D^2$ along the common part of the
boundary $\partial F\times D^2$ (respecting orientations of
course) to obtain a manifold (after smoothing the corner)
$$
W^{2n+2}:= N\bigcup_{\partial F\times D^2} (F\times D^2)~~~ {\rm
with}~~~ \partial W= F_{\psi}\cup (F\times S^1)\simeq {\cal
D}F_{\tp }$$ We use elementary obstruction theory to show that
this $W$ is stably parallelizable. Fix a frame field of the stable
tangent bundle of $N\subset W$. Obstructions to the extension of
this frame field over the whole manifold lie in the groups
$H^{q+1}(W,N;\pi_q(SO))\cong H^{q+1}(F,\partial F;\pi_q(SO))\cong
H_{2n-q-1}(F;\pi_q(SO))$. If $q\ne n-1$ or $q\ne 2n-1$ then
$H_{2n-q-1}(F;\zint)\cong 0$ (since $F$ has the homotopy type of
the wedge product of $n$-spheres). But if $q=n-1$ or $q=2n-1$ then
$\pi_q(SO)\cong 0$ because $n\equiv 3\pmod{8}$ and all
obstructions lie in the trivial groups. It implies that ${\cal
D}F_{\tp}$ is stable parallelizable and the Pontryagin class
$p_{\frac{n+1}{4}}({\cal D}F_{\tp})$ vanishes. Thus $\chi(\tp)=0$
(recall Lemma 2 of \cite{Kreck}) and hence $\chi_r(\psi)=0$.
\end{proof}

Now we can prove the following theorem of Stevens including the
cases when $n=3,7$ (cf. \cite{Stev}, Theorem 7).

\begin{cor}
Let $M_k$ be the k-fold branched cyclic cover of a
$(2n+1)$-manifold $M = F_{\vf}\cup_r (\partial F\times D^2)$ with
simple open book decomposition, where n is odd, $\ne 1$. Suppose
${\rm Var}(\vf^d)=0$, then $M_k$ and $M_{k+d}$ are (orientation
preserving) homeomorphic, while $M_k$ and $M_{d-k},~d>k$, are
orientation reversing homeomorphic.
\end{cor}
\begin{proof}
Since ${\rm Var}(\vf^d)=0$ and $M_d$ bounds a parallelizable
manifold (see Lemma 5 of \cite{Stev}) we have $\chi_r(\vf^d)=0$ by
the previous theorem. The exact sequence (\ref{kervar}) implies
that $\vf^d$ is isotopic to a diffeomorphism which belongs to the
image $\iota(\Theta_{2n+1})$ and therefore $F_{\vf^{d+k}}$ is
diffeomorphic to $F_{\vf^k}\# \Sigma'$ (cf. Lemma 1 of
\cite{Browd1}) for some $\Sigma'\in \Theta_{2n+1}$. In particular,
it means that $F_{\vf^{d+k}}$ is homeomorphic (via some
homeomorphism that preserves orientation) to $F_{\vf^{k}}$, and
hence $M_{d+k}$ is homeomorphic to $M_k$. To see the orientation
reversing case, notice that the mapping torus $F_g$ is
diffeomorphic to $F_{g^{-1}}$ via an orientation reversing
diffeomorphism induced, for instance, by the map $(x,t)\lmt
(g(x),1-t)$ from $F\times I$ to $F\times I$. This diffeomorphism
extends to an orientation reversing homeomorphism of the
corresponding open books $M$ and $M_{-1}$. Hence in our situation
$M_k=F_{\vf^{k}}\cup_r (\partial F\times D^2)$ is homeomorphic
(orient. revers.) to $F_{\vf^{-k}}\cup_r (\partial F \times D^2)$
which is homeomorphic (orient. pres.) to
$F_{\vf^{-k}}\#\Sigma'\cup_r (\partial F\times D^2)\simeq
F_{\vf^{d-k}}\cup_r (\partial F\times D^2)=M_{d-k}$.
\end{proof}

\noindent \underline{Remark}: If one defines $M_k :=
F_{\vf^k}\cup_r (\partial F\times D^2)$ for any $k\in \int$, then
the first statement that $M_k$ is homeomorphic to $M_{k+d}$
remains true, and the restriction $d>k$ in the second part can be
omitted.

To show diffeomorphism type periodicity we will basically use the
same argument plus the fact that the homotopy sphere $\Sigma'$
bounds a parallelizable manifold. We start with proving this fact.
Thus for $n\in\nat,~n\geq2$, we consider a diffeomorphism $h$ with
$[h]\in\K(V) $ such that our simple open book $M^{2n+1}=F_h\cup_r
(\partial F\times D^2)$ bounds a simply connected parallelizable
manifold $N$ and $\chi_r(h)=0$. In particular, we can assume that
$h\in {\rm Im}(\iota)$ is the identity except on a small closed
disk ${\cal D}^{2n}\hra {\rm int}(F)$ embedded into the interior
of $F$.

\begin{lem}
The natural inclusions $i_1:F\hra F_h$ and $i_2: F_h\hra M$ induce
isomorphisms $i_{1*}:H_n(F) \to H_n(F_h)$ and $i_{2*}:H_n(F_h) \to
H_n(M)$ respectively, and every $[z]\in H_n(M)$ can be represented
by an embedded sphere $S^n\hra M$ with trivial normal bundle
$\nu(S^n;M)$. In addition, $H_n(M)\cong H_{n+1}(M)$.
\end{lem}
\begin{proof}
That $i_{1*}$ is an isomorphism follows immediately from the Wang
exact sequence, and the other two isomorphisms follow from the
exact sequence of Stevens:
$$
0\lra H_{n+1}(M) \lra H_n(F,\partial F)\stackrel{\rm Var(h)}{\lra}
H_n(F)\lra H_n(M) \lra 0$$ which arises from the exact sequence of
the pair $(M,F)$ (see Proposition 1 of \cite{Stev}). Since the
normal bundle of any $S^n\hra M$ is stable and $M$ bounds a
parallelizable manifold, the bundle $\nu(S^n;M)$ must be trivial.
\end{proof}

Now we would like to kill $H_n(M)$ using surgery, and as a result
obtain a homotopy sphere $\Sigma_h\in \Theta_{2n+1}$ (we again
assume $n\geq 3$). We will show firstly that $\Sigma_h$ belongs to
$bP_{2n+2}$, and secondly that $\iota(\Sigma_h) =[h]$. Our
construction will follow Kreck's construction of the isomorphism
$\sigma: {\rm ker}(\chi) \lra \Theta_{2n+1}$ (see \cite{Kreck},
proof of Proposition 3).

For each generator $[z_i]\in H_n(F)$ we fix an embedding $\phi_i:
S^n_i\times d^{n+1}_i\hra F\times (0,1)\hra M$ disjoint from
${\cal D}^{2n}\times (0,1)$. Then we attach handles
$D^{n+1}_i\times d^{n+1}_i$ to the product $M\times I$ along these
embeddings into $M\times\{1\}$ to obtain a cobordism $W$ between
$M=M\times\{0\}$ and the homotopy sphere $\Sigma_h$ which is the
result of these $\phi_i$-surgeries on $M$. Furthermore, we can
choose the embeddings $\phi_i$ compatible with the framing of $M$
that comes from the framing of $N$ (see Lemma 6.2 of
\cite{KerMil}), and hence we get $W$ as a framed manifold. Taking
the union of $N$ and $W$ along $M$ we obtain a parallelizable
manifold with boundary $\Sigma_h$, and hence $\Sigma_h\in
bP_{2n+2}$. In the next lemma we show that $\Sigma_h$ is well
defined (depends only on the isotopy class of $h$) and that
$\iota(\Sigma_h)=[h]$, which implies $[h]\in \iota(bP_{2n+2})$
(cf. the properties of $W$ from \cite{Kreck}, pp. 655-656).

\begin{lem}~
\begin{enumerate}
\item The manifold $W$ is parallelizable, $n$-connected and
$\partial W=M\sqcup -\Sigma_h$. \item The embedding ${\cal L}:
M\times\{0\}\hra W$ induces an isomorphism\\ ${\cal L}_*:
H_{n+1}(M) \lra H_{n+1}(W)$ and all elements of $H_{n+1}(W)$ can
be represented by embedded spheres with trivial normal bundle.
\item If $W'$ is another n-connected manifold that also satisfies
property 2 and $\partial W'=M\sqcup -\Sigma$ for some
$\Sigma\in\Theta_{2n+1}$, then $\Sigma \simeq \Sigma_h$. \item
$\iota(\Sigma_h)=[h]$ where $\iota: \Theta_{2n+1}\hra \K(V) $.
\end{enumerate}
\end{lem}
\begin{proof}
$n$-connectivity follows immediately from $(n-1)$-connectivity of
$M$ and the Mayer-Vietoris exact sequence of the union
$$
W:=M\times [0,1] \bigcup (\sqcup_i D^{n+1}_i\times d^{n+1}_i)
$$ Hence, by the Hurewicz theorem we can represent every $[z]\in
H_{n+1}(W)$ by an embedded sphere $S^{n+1}\hra W$. The same exact
sequence implies that the embedding ${\cal L}: M=M\times\{0\} \hra
W$ induces isomorphism between $H_{n+1}(M)$ and $H_{n+1}(W)$ and
to finish part 2 we just need to show that $\nu(S^{n+1};W)$ is
trivial. To see this, notice that every $[z]\in H_{n+1}(M)\cong
H_{n+1}(W)$ can be represented by an embedded $S^n\times S^1\hra
M$ (see the Lemma above) with trivial normal bundle in $M\times
(0,1)$. Since we can take the diffeomorphism $h$ to be the
identity on $\nu(S^n;F)$, it is not hard to see that there is an
embedded $S^{n+1}\hra M$ which is cobordant to $S^n\times S^1$ in
$M$ and hence framed cobordant in $M\times(0,1)\hra W$. This
proves 2. Suppose now we have $W'$ that satisfies property 3. Take
a union of $W$ and $W'$ along $M$ to obtain a n-connected
cobordism $\cal W$  between $\Sigma$ and $\Sigma_h$. We will show
that we can make it into an $h$-cobordism. Mayer-Vietoris exact
sequence implies that $H_{n+1}({\cal W})\cong H_{n+1}(W)\oplus
H_n(M)$, and since $H_{n+1}(W)\cong H_{n+1}(W,M)$ (Poincar\'e
duality plus Universal coefficient theorem) and also
$H_{n+1}(W,M)\cong H_n(M)$ we see that $H_{n+1}({\cal W})$ has the
direct summand $H_{n+1}(W)$ with
the properties: \\
i) ${\rm dim}(H_{n+1}(W))=\frac{1}{2}{\rm dim}(H_{n+1}({\cal W}))$\\
ii) every homology class of $H_{n+1}(W)$ can be represented by an
embedded sphere $S^{n+1}\hra W$ with trivial normal bundle.\\
iii) For all $[z_1],~[z_2]\in H_{n+1}(W)$ the intersection number
$z_1\circ z_2$ vanishes (this follows from 2. of this lemma).
Therefore we can use surgery to kill $H_{n+1}({\cal W})$ and
obtain an $h$-cobordism between $\Sigma$ and $\Sigma_h$. The last
property follows from \cite{Browd1}, Lemma 1. Indeed, if
$\Sigma'=D^{2n+1}_1\cup_h D^{2n+1}_2$ then by our definition
$\iota(\Sigma')=[h]\in \K(V) $ and hence $M\simeq\partial(F\times
D^2)\# \Sigma'$. Since $\partial(F\times D^2)$ is framed cobordant
to the standard $(2n+1)$-sphere, we can use this cobordism (namely
$F\times D^2$ minus an embedded disk $D^{2n+2}$) to produce a
$n$-connected cobordism between $M$ and $\Sigma'$ that will
satisfy property 2. As we have just seen above this means that
$\Sigma'\simeq \Sigma_h$, i.e. $\iota(\Sigma_h) = [h]$.
\end{proof}

Let us denote the signature of a parallelizable manifold $N_k$
with boundary $\partial N_k= M_k = F_{\vf^k}\cup_r (\partial
F\times D^2)$ by $\sigma_k$, and the generator of $bP_{2n+2}$ by
$\Sigma$.

\begin{cor}(cf. \cite{Stev}, Proposition 8) Let $M^{2n+1}=
F_{\vf}\cup_r (\partial F\times D^2)$ be the manifold with simple
open book decomposition where $n$ is odd, $\ne 1$ and $M_k$ be the
$k$-fold branched cyclic cover of $M$. If ${\rm Var}(\vf^d)=0$
then $M_{k+d}$ is diffeomorphic to
$(\frac{\sigma_d}{8}\cdot\Sigma)\# M_k$.
\end{cor}
\begin{proof}
We have just seen above that $M_{k+d}\simeq \Sigma'\# M_k$ with
$\Sigma'=m\cdot\Sigma\in bP_{2n+2}$ for some $m\in\nat$. Since
$M_d=F_{\vf^d}\cup_r (\partial F\times D^2)\simeq \partial(F\times
D^2)\# m\Sigma$ and $m\Sigma$ bounds a parallelizable manifold,
say $W_m$, with signature $\sigma(W_m)=8m$ and $\partial(F\times
D^2)$ bounds $F\times D^2$ (which is also parallelizable) with
signature zero, the connected sum of $W_m$ and $F\times D^2$ along
the boundary (cf.~\S 2~of \cite{KerMil}) will give us a
parallelizable manifold $N_d := W_m\# F\times D^2$ with boundary
$\partial N_d = M_d$ and signature
$\sigma(N_d)=\sigma(W_m)+\sigma(F\times D^2)=8m+0$. Thus
$m=\frac{\sigma(N_d)}{8}\equiv\frac{\sigma_d}{8}\bmod ({\rm
order~of~}bP_{2n+2})$ and the corollary follows.
\end{proof}

When $n=even$, the periodicity of $M_k$ is more complicated.
Consider the link of the singularity $z^2_0+z^2_1+\ldots +z^2_n=0$
with $n=2m$ and denote the $(4m+1)$-dimensional Kervaire sphere by
$\Sigma$ and the tangent $S^n$-sphere bundle to $S^{n+1}$ by $T$.
Then $M_{k+8}$ is diffeomorphic to $M_k$ and the diffeomorphism
types are listed in the table (see \cite{DK}, Proposition 6.1)

\parskip=5mm

\begin{tabular}{|c|c|c|c|c|c|c|c|}
\hline $M_1\simeq M_7$ & $M_2$ & $M_3\simeq M_5$ & $M_4$ & $M_6$ &
$M_8$ \\ \hline $S^{2n+1}$ & ~$T$~ & ~$\Sigma$~ & ~$(S^n\times
S^{n+1})\#\Sigma$~ & ~ $T\#\Sigma$~ &
 ~ $S^n\times S^{n+1}$\\
\hline
\end{tabular}

\noindent The following result is due to Stevens (\cite{Stev},
Theorem 9).

\begin{cor}
If for branched cyclic covers $M_k$ of a $(2n+1)$-manifold $M$
with simple open book decomposition ${\rm Var}(\vf^d)=0$, then
$M_k$ and $M_{k+2d}$ are homeomorphic and $M_k$ and $M_{k+4d}$ are
diffeomorphic. Moreover, if $n=2~or~6$, then $M_k$ and $M_{k+d}$
are diffeomorphic.
\end{cor}
\begin{proof}
When $n=2$ the mapping class group is trivial and $[\vf^d]=Id$. If
$n=6$ then $G\cong 0$ and $bP_{14}\cong 0$ (see \cite{KerMil},
Lemma 7.2) which implies that $[\vf^d]=Id$. For the other even $n$
we know that the group $G$ is isomorphic either to $\zint_2$ or
$\zint_2\oplus \zint_2$ and hence $\chi_r(\vf^d)$ has order two.
Therefore $\vf^{2d}\in bP_{4m+2}$, i.e. $M_k$ is homeomorphic to
$M_{k+2d}$, and since the group $bP_{4m+2}$ is either trivial or
$\zint_2$ (see \cite{KerMil}), $\vf^{4d}$ must be pseudo-isotopic
to the identity.
\end{proof}

\noindent {\tt Example 2:} (The authors are indebted to the
referee for suggesting this example.) Consider again the
singularity $z^2_0+z^2_1+\ldots +z^2_n=0$ with $n=2m$. Assume in
addition that $m\ne 0$ (mod 4) and that the Kervaire sphere
$\Sigma \in bP_{4m+2}$ is exotic, e.g. when $4m+2\ne 2^l-2$ (see
\cite{Browder}). Here the Milnor fiber $F$, is the tangent disc
bundle to the sphere $S^{2m}$ and hence $\DblF\simeq S^{2m}\times
S^{2m}$. It is also well known that the geometric monodromy $\vf$
of this singularity satisfies the properties: $\vf_*=-Id$, ${\rm
Var}(\vf^2)=0$ and ${\rm Var}(\vf)\ne 0$ (cf. \cite{Looijenga},
Chapter 3). Since $M_0$ is not diffeomorphic to $M_2$ and $M_1$ is
not diffeomorphic to $M_5$, $\chi_r([\vf^2])$ will be a generator
of $\Rhom\cong\int_2$ and $[\vf^4]$ will be a generator of
$bP_{4m+2}\cong\int_2$. Since $\Theta_{4m+1}\cong bP_{4m+2}\oplus
{\rm Coker}(J_{4m+1})$ (cf. \cite{Brumfiel}) we see that in this
case $\K(V) \cong\int_4\oplus {\rm Coker}(J_{4m+1})$ and the exact
sequence (\ref{kervar}) doesn't split.

\subsection{Periodicity in dimension 3}

It is known that if the dimension of the open book $M^{2n+1}$ is
three, then there is homological periodicity (see references in
\cite{DK}) but there is no topological one. For the sake of
completeness we illustrate this with the following classical
example (cf. \cite{Rolfsen}, Chapter 10. D.). Let
$f(z_0,z_1)=z_0^2+z_1^3$ be the complex polynomial which defines
the curve $V=\{f(z_0,z_1)=0\}$ in $\cn^2$ with the cusp at the
origin. The corresponding Milnor fibration has monodromy $\vf$ of
order six, the boundary of the fiber $F$ is the trefoil knot $K$
and $\Va (\vf^6)=0$. This fibration gives the open book structure
to the standard 3-sphere $S^3=M_1=F_{\vf}\cup (K\times D^2)$. We
show that $M_7\ne M_1$ and $M_6\ne M_0=(F\times S^1)\cup (K\times
D^2)$.

Let us first compare $\pi_1(M_0)$ with $\pi_1(M_6)$. The theorem
of Seifert and Van Kampen entails that $\pi_1(M_0)\cong \pi_1(F)$
which is the free group on two generators. As for $M_6$ one can
easily find using the Reidemeister-Schreier theorem a presentation
for $\pi_1(F_{\vf^6})$ and then show that $\pi_1(M_6)$ admits the
following presentation:
$$
\left\langle
Z_1,Z_2,\ldots,Z_6~|~Z_1=Z_6Z_2,~\ldots,~Z_j=Z_{j-1}Z_{j+1},~
\ldots,~ Z_6=Z_5Z_1\right\rangle
$$
It takes a bit more effort to show that this group is isomorphic
to the group of upper unitriangular $3\times 3$ matrices with
integer coefficients (cf. \cite{Mil3}, \S8)
$$
{\cal H}\cong \{
\begin{pmatrix}
1 & a & c\\
0 & 1 & b\\
0 & 0 & 1
\end{pmatrix}|a,b,c\in\int \}
$$

Suppose now that $M_7$ were homeomorphic to the sphere. Then we
could take the union of $N_7$ and $D^4$ (recall that $N_7$ is the
cyclic covering of $D^4$ branched along the fiber $(F,K)\hra
(D^4,S^3)$ where $F\cap S^3=K$):
$$
W^4:=N_7\bigcup_{S^3} D^4
$$
Since $N_7$ is parallelizable (see \cite{Cappell}, Theorem 5 or
\cite{Kauf-book}, Chapter XII), $W^4$ would be a closed
spin-manifold. Hence its signature $\sigma(W^4)=\sigma(N_7)$ must
be a multiple of 16 by the theorem of Rokhlin \cite{Rokh}. But
$\sigma(N_7)=-8$ as one can find using the Seifert pairing on
$H_1(F)$ (cf. \cite{Cappell}; \cite{Kauf-book}), and hence $M_7\ne
M_1$. Actually much more is known. Milnor in \cite{Mil3} proved
that $\pi_1(M_r)$ is isomorphic to the commutator subgroup
$[\Gamma,\Gamma]$ of the centrally extended triangle group
$\Gamma$ which has a presentation
$$
\Gamma\cong \langle
\gamma_1,~\gamma_2,~\gamma_3~|~\gamma_1^2=\gamma_2^3=
\gamma_3^r=\gamma_1\cdot \gamma_2\cdot \gamma_3 \rangle
$$
This group $\Gamma$ is infinite when $r\geq 6$ (see \cite{Mil3},
\S2,3) and hence $[\Gamma,\Gamma]$, that has index $r-6$, is
infinite too. In particular, none of the cyclic coverings of $S^3$
branched along the trefoil knot can be simply connected.


\parskip=1mm

\noindent University of Illinois at Chicago \\
Department of Mathematics, Statistics and Computer Science\\
851 S.Morgan st. Chicago, IL 60607

\noindent {\small E-mail address:~kauffman@uic.edu}

~

\noindent International University Bremen\\
School of Engineering \& Science\\
P.O. Box 750 561\\
28725 Bremen, Germany

\noindent {\small E-mail address:~n.krylov@iu-bremen.de}

\end{document}